\documentclass[preprint,12pt,3p]{elsarticle}

\usepackage{amssymb}
\usepackage{amsmath,amsthm}
\usepackage{mathrsfs}
\usepackage{hyperref}
\usepackage{cleveref}
\usepackage{needspace}
\usepackage[all,cmtip]{xy}
\usepackage{epsf, graphicx}
\usepackage{latexsym,amsfonts,amsbsy,amssymb}
\usepackage{geometry}
\usepackage{titletoc}
\usepackage{color,xcolor}
\usepackage{rotating}
\usepackage{algorithm}
\usepackage{algorithmic}
\usepackage{amscd}
\makeatletter
\renewcommand\theequation{\thesection.\arabic{equation}}
\@addtoreset{equation}{section} \makeatother
\newtheorem{theorem}{Theorem}[section]
\newtheorem{lemma}[theorem]{Lemma}
\newtheorem{notation}[theorem]{Notation}
\newtheorem{conjecture}[theorem]{Conjecture}
\newtheorem{void}[theorem]{}
\newtheorem{proposition}[theorem]{Proposition}%[section]
\def\IBr{{\rm IBr}}

\def\Ind{{\rm Ind}}
\def\Res{{\rm Res}}

\def\br{{\rm br}}

\def\Aut{{\rm Aut}}

\def\rad{{\rm rad}}

\def\End{{\rm End}}
\def\Gal{{\rm Gal}}

\def\W{\mathcal{W}}
\def\H{\mathcal{H}}
\def\F{\mathbb{F}}

\def\hd{{\rm hd}}

\makeatletter
\def\ps@pprintTitle{%
\let\@oddhead\@empty
\let\@evenhead\@empty
\def\@oddfoot{\reset@font\hfil\thepage\hfil}
\let\@evenfoot\@oddfoot
}
\newcommand{\dedication}[1]{\g@addto@macro\elsauthors{\par\medskip{\normalsize\itshape #1}\par}}
\makeatother

\begin{document}

\begin{frontmatter}
	
	\title{Alperin weight bijections and dimensions of simple modules}

	\author{Xin Huang}

	%\dedication{}
	%\fntext[label1]{The author acknowledges the support.}
	
%	\address{School of Mathematics and Statistics, Central China Normal University, Wuhan 430079, China}
%	\ead{xinhuang@mails.ccnu.edu.cn}

\author{Hanxiao Li}
	
	\begin{abstract}
We propose a new refinement of the (blockwise) Alperin weight conjecture taking dimensions of simple modules (or degrees of irreducible Brauer characters) into consideration. This is similar to Giannelli's refinement of the (Alperin--)McKay conjecture. We prove this new refinement for blocks with a cyclic or Klein four defect group, for inertial blocks (which include the Brou\'{e}--Puig nilpotent blocks), and for $p$-blocks of $p$-solvable groups.
\end{abstract}
	
	
\end{frontmatter}

% \linenumbers

\section{Introduction}\label{s1}
Throughout this paper $p$ is a prime number and $k$ is an algebraic closure of $\F_p$. Let $G$ be a finite group. Denote by $\IBr(G)$ the set of isomorphism classes of simple $kG$-modules, and by $\W(G)$ the set of isomorphism classes of projective simple $kG$-modules. Alperin's weight conjecture predicts that there exists a bijection between
$$\IBr(G)~~~{\rm and}~~~\bigsqcup_{(P)}\W(N_G(P)/P),$$
where the disjoint union is taken over a set of representatives of the conjugacy classes of $p$-subgroups of $G$. 
There is also a ``blockwise" version of this conjecture which will be disscussed in detail later. This conjecture has been one of the central problems in the representation theory of finite groups for nearly 40 years.

Over the years, several refinements of the (blockwise) Alperin weight conjecture have been proposed, aiming to provide a deeper understanding of the relationship between $\IBr(G)$ and $\bigsqcup_{(P)}\W(N_G(P)/P)$. Notably, Navarro \cite{Nav04} suggested a refinement taking Galois automorphisms (or equivalently, fields of values) into account. The inductive (blockwise) Alperin weight condition \cite{NT, Spa} can also be viewed a strong form of the (blockwise) Alperin weight conjecture taking more information, including the action of group automorphisms, into consideration. There is an even stronger version of the (blockwise) Alperin weight conjecture, which considers both actions of Galois automorphisms and group automorphisms, and which has been studied in \cite{DH,DHZ,FFZ1,FFZ2,H26}. These refined conjectures not only generalize Alperin's original statement but also suggest a rich and intricate interplay between local and global representation-theoretic invariants.

Inspired by a recent refinement of the McKay conjecture proposed by Giannelli \cite{G}, we propose a further strengthening of the Alperin weight conjecture:

\begin{conjecture}\label{conj: non-block version}
Let $G$ be a finite group and let $s_p(G)$ be the set of all $p$-subgroups of $G$. Then there exists a bijection
$$f:\IBr(G)\to\bigsqcup_{P\in [s_p(G)/\sim_G]}\W(N_G(P)/P),$$
such that for any isomorphism class $[V]$ of a simple $kG$-module $V$, we have $\dim_k(f([V]))\leq \dim_k([V])$. Here $[s_p(G)/\sim_G]$ denotes a set of representatives of $s_p(G)$ under $G$-conjugation, and $\dim_k([V])$ is defined to be $\dim_k(V)$.
\end{conjecture}

It is natural to believe that there should be a block version of Conjecture \ref{conj: non-block version}. To state it, we need some terminology and notation. Let $G$ be a finite group. For any $p$-subgroup $P$ of $G$, let $(kG)^P=\{a\in kG\mid uau^{-1}=a~{\rm for~all}~u\in P\}$. The {\it $P$-Brauer homomorphism} $\br_P^{kG}:(kG)^P\to kC_G(P)$ is the map sending any $\sum_{g\in G}\alpha_gg \in (kG)^P$ to $\sum_{g\in C_G(P)}\alpha_gg$, where $\alpha_g\in k$. It is well-known that the $P$-Brauer homomorphism is a homomorphism of $k$-algebras, and it restricts to a $k$-algebra homomorphism $Z(kG)\to (kC_G(P))^{N_G(P)}=\{a\in kC_G(P)\mid hah^{-1}=a~{\rm for~all}~h\in N_G(P)\}$. A {\it block} $b$ of $kG$ is a primitive central idempotent of $kG$, and $kGb$ is called a {\it block algebra} of $kG$. For any central idempotent $b$ of $G$, we denote by $\IBr(G,b)$ the set of isomorphism classes of simple $kGb$-modules, and by $\W(G,b)$ the set of isomorphism classes of projective simple $kGb$-modules. The following is the block-theoretic version of Conjecture \ref{conj: non-block version}:

\begin{conjecture}\label{conj:block version}
Let $G$ be a finite group and $b$ a central idempotent of $kG$. Denote by $\bar{b}_P$ the image of $b_P:=\br_P^{kG}(b)$ in $kN_G(P)/P$, which is a central idempotent of $kN_G(P)/P$. Then there exists a bijection
$$f:\IBr(G,b)\to\bigsqcup_{P\in [s_p(G)/\sim_G]}\W(N_G(P)/P,\bar{b}_P),$$
such that for any isomorphism class $[V]$ of a simple $kGb$-module $V$, we have $\dim_k(f([V]))\leq \dim_k([V])$. 
\end{conjecture}

If we take $b=1$, then Conjecture \ref{conj:block version} returns to Conjecture \ref{conj: non-block version}. It is clear that in order to prove Conjecture \ref{conj:block version}, it suffices to verify it for all block idempotents $b$ of $kG$. To provide evidence of Conjecture \ref{conj:block version}, we prove the following three results, which include several infinite families of blocks.

\begin{theorem}\label{theo:cyclic}
Let $G$ be a finite group and $b$ a block of $kG$ with a defect group $D$. Let $\H= \Gal(k/\F_p)\times \Aut(G)$. Assume that $D$ is a cyclic group or a Klein four group. Then there exists a bijection
$$f:\IBr(G,b)\to\bigsqcup_{P\in [s_p(G)/\sim_G]}\W(N_G(P)/P,\bar{b}_P),$$
satisfying the following two conditions:
\begin{enumerate}[{\rm (i)}]
	\item For any $[V]\in \IBr(G,b)$, we have $\dim_k(f([V]))\leq \dim_k([V])$.
	\item $f$ commutes with the action of $\H_b$, the stabilizer of $b$ in $\H$.
\end{enumerate}
\end{theorem}

Theorem \ref{theo:cyclic} provides evidence for an even stronger refinement of the blockwise Alperin weight conjecture which considers dimensions, Galois automorphisms and group automorphisms. The refinement which considers only Galois automorphisms was proved for blocks with a cyclic or Klein four defect group; see \cite[Corollary 3 (i)]{Huang2024}. So the above theorem generalizes \cite[Corollary 3 (i)]{Huang2024}. For simplicity, in the remaining cases, we mainly provide evidence for the dimension refinement Conjecture \ref{conj:block version}, and we don't pursue much generality. 

Following Puig \cite[1.4]{Puig11}, a block algebra $kGb$ with a defect group $D$ is called {\it inertial} if it is Morita equivalent to its Brauer correspondent $kN_G(D)c$ via a bimodule with an endopermutation source. The nilpotent blocks defined by Brou\'{e} and Puig \cite[Definition 1.1]{BP} is a special class of inertial blocks; see \cite[1.5]{Puig11}.

\begin{theorem}\label{theo:inertial blocks}
	Let $G$ be a finite group and $b$ an inertial block of $kG$. Then there exists a bijection
	$$f:\IBr(G,b)\to\bigsqcup_{P\in [s_p(G)/\sim_G]}\W(N_G(P)/P,\bar{b}_P)$$
such that for any $[V]\in \IBr(G,b)$, we have $\dim_k(f([V]))\leq \dim_k([V])$.
\end{theorem}

Let $G$ be a finite group and $[M]\in\IBr(G)$.
We define the {\it codegree}
of $[M]$ by
${\rm codeg}([M])=|G|/\dim_k(M)$
and denote by ${\rm codeg}([M])_{p'}$
the $p'$-part of ${\rm codeg}([M])$.

\begin{theorem}\label{theo:p-solvable groups}
	Let $G$ be a finite $p$-solvable group and $b$ a block of $kG$. Then there exists a bijection
	$$f:\IBr(G,b)\to\bigsqcup_{P\in [s_p(G)/\sim_G]}\W(N_G(P)/P,\bar{b}_P),$$
    satisfying the following conditions:
	\begin{enumerate}[{\rm (i)}]
		\item For any $[V]\in \IBr(G,b)$, we have $\dim_k(f([V]))\leq \dim_k([V])$.
        \item For any $[V]\in \IBr(G,b)$, we have ${\rm codeg}([V])_{p'}^2 \equiv {\rm codeg} (f([V]))_{p'}^2 \pmod p$.
		\item $f$ commutes with the action of $\Gal(k/\F_p)_b$, the stabilizer of $b$ in $\Gal(k/\F_p)$.
	\end{enumerate}
\end{theorem}

We remark that the codegree condition in Theorem \ref{theo:p-solvable groups} (ii) does not hold in general for arbitrary finite groups; see the first paragraph of \cite[Introduction]{Turull14} for details. 

In Section \ref{s2} we review the actions of Galois automorphisms and group automorphisms on modules and fix some related notation. Theorems \ref{theo:cyclic}, \ref{theo:inertial blocks} and \ref{theo:p-solvable groups} will be proved in Sections \ref{s3}, \ref{s4} and \ref{s5}, respectively.

\section{The actions of Galois automorphisms and group automorphisms}\label{s2}
Throughout this section $G$ is a finite group and $b$ is a central idempotent of $kG$. The motivation of this section is to explain the actions of $\Gal(k/\F_p)$ and $\Aut(G)$ on the sets $$\IBr(G)~~~{\rm and}~~~\bigsqcup_{P\in [s_p(G)/\sim_G]}\W(N_G(P)/P),$$ and also to fix some notation. Some of them were explained in \cite{H26} in a more general context. We formulate for conventient reference in this section that the actions of $\Gal(k/\F_p)$ and $\Aut(G)$ are compatible with the Green correspondence.

\begin{void}
{\rm \textbf{The action of Galois automorphisms.} (i) Let $V$ be a $kG$-module and $\sigma\in\Gal(k/\F_p)$. We denote by ${}^\sigma V$ the $kG$-module, which is equal to $V$ as a module over the subring $\F_pG$ of $kG$, such that $x$ acts on ${}^\sigma V$ as $\sigma^{-1}(x)$ for all $x\in k$. Note that if $f:V_1\to V_2$ is a homomorphism (resp. isomorphism) of $kG$-modules, then $f:{}^\sigma V_1\to {}^\sigma V_2$ is also a homomorphism (resp. isomorphism) of $kG$-modules. Hence $\sigma$ defines an auto-functor of the category of $kG$-modules. Moreover, for the isomorphism class $[V]$ of $V$, we can define ${}^\sigma[V]$ to be $[{}^\sigma V]$. It is clear that ${}^\sigma V$ is simple if and only if $V$ is simple. Hence the group $\Gal(k/\F_p)$ acts on the set $\IBr(G)$. One checks that the map $kG\to {}^\sigma(kG)$ sending any $\sum_{g\in G}\alpha_gg$ to $\sum_{g\in G} \sigma^{-1}(\alpha_g)g$ is an isomorphism of $kG$-modules. Hence $V$ is projective if and only if ${}^\sigma V$ is projective. Therefore, the group $\Gal(k/\F_p)$ also acts on the set $\W(G)$. Hence it acts on $\bigsqcup_{P\in [s_p(G)/\sim_G]}\W(N_G(P)/P)$ by acting on each $\W(N_G(P)/P)$.

(ii) The field automorphism $\sigma$ induces an $\F_p$-algebra automorphism of $kG$, sending any $\sum_{g\in G}\alpha_gg$ to $\sum_{g\in G}\sigma(\alpha_g)g$; we abusively use the same letter $\sigma$ to it. If $V$ is a $kG$-module and $kG\to \End_k(V)$ is the structural homomorphism of $V$, then the structural homomorphism of ${}^\sigma V$ is $kG\xrightarrow{\sigma^{-1}} kG\to \End_k(V)$. If $V$ is a $kGb$-module, then ${}^\sigma V$ is a $kG\sigma(b)$-module. Hence $\sigma$ induces a bijections which we still denote them by the same letter:
$$\sigma:\IBr(G,b)\to \IBr(G,\sigma(b)),~~~[V]\mapsto [{}^\sigma V].$$
$$\sigma:\W(G,b)\to \W(G,\sigma(b)),~~~[U]\mapsto [{}^\sigma U].$$
Let $P$ be a $p$-subgroup of $G$. It is straightforward to check the following commutative diagram:
$$\xymatrix{ (kG)^P\ar[d]_{\sigma}\ar[rr]^{\br_P^{kG}} & & kC_G(P)\ar[d]^{\sigma}\ar@{^{(}->}[rr] &  & kN_G(P)\ar[d]^{\sigma} \ar@{->>}[rr]  &  & kN_G(P)/P \ar[d]^{\sigma}\\
	(kG)^P\ar[rr]_{\br_P^{kG}} &  & kC_G(P)\ar@{^{(}->}[rr]  &  &   kN_G(P) \ar@{->>}[rr] &   & kN_G(P)/P
}$$
Recalled that we denoted by $b_P$ the element $\br_P^{kG}(b)$ and by $\bar{b}_P$ the image of $b_P$ in $kN_G(P)/P$. By the above commutative diagram we have $\sigma(b_P)=\sigma(b)_P$ and $\sigma(\bar{b}_P)=\overline{\sigma(b)}_P$.
Hence $\sigma$ induces a bijection which we still denote it by the same letter:
$$\sigma:\W(N_G(P)/P,\bar{b}_P)\to \W(N_G(P)/P,\overline{\sigma(b)}_P),~~~[U]\mapsto [{}^\sigma U].$$
}
\end{void}

\begin{void}\label{void:The action of group automorphisms}
{\rm \textbf{The action of group automorphisms.} (i) Let $V$ be a $kG$-module. Let $H$ be another finite group which is isomorphic to $G$, and let $\varphi:G\to H$ be a group isomorphism. We denote by ${}^\varphi V$ the $kH$-module, which is equal to $V$ as a $k$-module, such that $h$ acts on ${}^\varphi V$ as $\varphi^{-1}(h)$ for all $h\in H$. Note that if $f:V_1\to V_2$ is a homomorphism (resp. isomorphism) of $kG$-modules, then $f:{}^\varphi V_1\to {}^\varphi V_2$ is a homomorphism (resp. isomorphism) of $kH$-modules. Hence $\varphi$ defines a functor from the category of $kG$-modules to the category of $kH$-modules. Moreover, for the isomorphism class $[V]$ of $V$, we can define ${}^\varphi[V]$ to be $[{}^\varphi V]$. It is clear that ${}^\varphi V$ is a simple $kH$-module if and only if $V$ is a simple $kG$-module. Hence the group $\Aut(G)$ acts on the set $\IBr(G)$. One checks that the map $kG\to {}^\varphi(kG)$ sending any $\sum_{g\in G}\alpha_gg$ to $\sum_{g\in G} \alpha_g\varphi^{-1}(g)$ is an isomorphism of $kG$-modules. Hence $V$ is projective if and only if ${}^\varphi V$ is projective. Therefore, the group $\Aut(G)$ also acts on the set $\W(G)$.

(ii) The group isomorphism $\varphi$ induces a $k$-algebra isomorphism of $kG\to kH$, sending any $\sum_{g\in G}\alpha_gg$ to $\sum_{g\in G}\alpha_g\varphi(g)$; we abusively use the same letter $\sigma$ to it. If $V$ is a $kG$-module and $kG\to \End_k(V)$ is the structural homomorphism of $V$, then the structural homomorphism of ${}^\varphi V$ is $kH\xrightarrow{\varphi^{-1}} kG\to \End_k(V)$. If $V$ is a $kGb$-module, then ${}^\varphi V$ is a $kH\varphi(b)$-module. Hence $\varphi$ induces a bijections which we still denote them by the same letter:
$$\varphi:\IBr(G,b)\to \IBr(H,\varphi(b)),~~~[V]\mapsto [{}^\varphi V].$$
$$\varphi:\W(G,b)\to \W(H,\varphi(b)),~~~[U]\mapsto [{}^\varphi U].$$

(iii) Now let $\varphi\in \Aut(G)$ and $P$ a $p$-subgroup of $G$. The $\varphi$ restricts to group isomorphisms $C_G(P)\to C_G(\varphi(P))$, $N_G(P)\to N_G(\varphi(P))$ and $N_G(P)/P\to N_G(\varphi(P))/\varphi(P)$, which we denote them by the same letter $\varphi$. It is straightforward to check the following commutative diagram:
$$\xymatrix{ (kG)^P\ar[d]_{\varphi}\ar[rr]^{\br_P^{kG}~~~} & & kC_G(P)\ar[d]^{\varphi}\ar@{^{(}->}[rr] &  & kN_G(P)\ar[d]^{\varphi} \ar@{->>}[rr]  &  & kN_G(P)/P \ar[d]^{\varphi}\\
	(kG)^{\varphi(P)}\ar[rr]_{\br_{\varphi(P)}^{kG}~~~} &  & kC_G(\varphi(P))\ar@{^{(}->}[rr]  &  &   kN_G(\varphi(P)) \ar@{->>}[rr] &   & kN_G(\varphi(P))/\varphi(P)
}$$
By the above commutative diagram we have $\varphi(b_P)=\varphi(b)_{\varphi(P)}$ and $\varphi(\bar{b}_P)=\overline{\varphi(b)}_{\varphi(P)}$. Then by \ref{void:The action of group automorphisms} (ii), $\varphi$ induces a bijection which we still denote it by the same letter:
$$\varphi:\W(N_G(P)/P,\bar{b}_P)\to \W(N_G(\varphi(P))/\varphi(P),\overline{\varphi(b)}_{\varphi(P)}),~~~[U]\mapsto [{}^\varphi U].$$
	
(iv) For any $g\in G$, we denote by $c_g$ the inner automorphism of $G$ sending any $x\in G$ to ${}^gx=gxg^{-1}$. If $H$ is a subgroup of $G$, then we use the same symbol $c_g$ to denote the isomorphism $H\to {}^gH=gHg^{-1}$ sending any $h\in H$ to ${}^gh$. 

(v) For any $\varphi\in \Aut(G)$ and any $$[U]\in \bigsqcup_{P\in [s_p(G)/\sim_G]}\W(N_G(P)/P),$$ there is a unique $P\in [s_p(G)/\sim_G]$ such that $[U]\in \W(N_G(P)/P)$. There is also a unique $Q\in [s_p(G)/\sim_G]$ such that $\varphi(P)={}^gQ$ for some $g\in G$. Equivalently, $\varphi(P)=c_g(Q)$. Then we define
$${}^\varphi[U]:=(c_g^{-1}\circ \varphi)([U])=[{}^{(c_g^{-1}\circ \varphi)}U]\in \W(N_G(Q)/Q).$$
Note that $[{}^\varphi U]$ is independent of the choice of $g$. Indeed, if $\varphi(P)={}^hQ$ for another $h\in G$, then $t:=h^{-1}g\in N_G(Q)$ which induces an identity map $c_t$ on $\W(N_G(Q)/Q)$. Now we have explicitly described the action of $\Aut(G)$ on 
$$\bigsqcup_{P\in [s_p(G)/\sim_G]}\W(N_G(P)/P).$$
The resulting $\Aut(G)$-set is independent, up to isomorphism, of the choice of the representatives $[s_p(G)/\sim_G]$; see Proposition \ref{prop:independence of representatives}.

}
\end{void}

\begin{void}
{\rm \textbf{The mixed action of Galois and group automorphisms.} Let $H$ be another finite group which is isomorphic to $G$. For any $\sigma\in \Gal(k/\F_p)$ and any group isomorphism $\varphi:G\to H$, we have the following commutative diagram:
	$$\xymatrix{kG\ar[rr]^{\sigma}\ar[d]_{\varphi} &  & kG\ar[d]^{\varphi}\\
		kH \ar[rr]_{\sigma}  &   &  kH 
	}$$
Hence $\sigma(\varphi(b))=\varphi(\sigma(b))$ and we have the commutative diagram
	$$\xymatrix{\IBr(G,b)\ar[rr]^{\sigma}\ar[d]_{\varphi} &  & \IBr(G,\sigma(b))\ar[d]^{\varphi}\\
 \IBr(H,\varphi(b)) \ar[rr]_{\sigma}  &   &  \IBr(H,\sigma(\varphi(b))) 
}$$
Moreover, for any $p$-subgroup $P$ of $G$, we have the commutative diagram	
$$\xymatrix{\W(N_G(P)/P,\bar{b}_P)\ar[rr]^{\sigma}\ar[d]_{\varphi} &  & \W(N_G(P)/P,\overline{\sigma(b)}_P)\ar[d]^{\varphi}\\
	\W(N_H(\varphi(P))/\varphi(P),\overline{\varphi(b)}_{\varphi(P)}) \ar[rr]_{\sigma}  &   &  \W(N_H(\varphi(P))/\varphi(P),\overline{\sigma(\varphi(b))}_{\varphi(P)}) 
}$$
Therefore, the direct product $\H:=\Gal(k/\F_p)\times \Aut(G)$ acts on both sets 
$$\IBr(G)~~~{\rm and}~~~\bigsqcup_{P\in [s_p(G)/\sim_G]}\W(N_G(P)/P),$$
and this action restricts to an action of $\H_b$ on the sets 
$$\IBr(G,b)~~~{\rm and}~~~\bigsqcup_{P\in [s_p(G)/\sim_G]}\W(N_G(P)/P,\bar{b}_P).$$

}
\end{void}

\begin{notation}[The Green correspondence and the head of modules]\label{notation: green corr and head}
	{\rm Let $G$ be a finite group and $P$ a $p$-subgroup of $G$. Let $V$ be an indecomposable $kG$-module with a vertex $P$, we denote by $f_P(V)$ a Green correspondent of $V$; namely, $f_P(V)$ is an indecomposable direct summand of $\Res_{N_G(P)}^G(V)$ with a vertex $P$. We define $f_P([V])$ to be $[f_P(V)]$. Let $U$ be an indecomposable $kN_G(P)$-module with a vertex $P$, we denote by $g_P(U)$ a Green correspondent of $U$; that is, $g_P(U)$ is an indecomposable direct summand of $\Ind_{N_G(P)}^G(U)$ with a vertex $P$. We define $g_P([U])$ to be $[g_P(U)]$. Let $M$ be a finite dimensional $kG$-module. We write $\hd(M):= M / \rad(M)$ (where $\rad(M)$ is the radical of $M$), which is called the {\it head} of $M$; we define $\hd([M])$ to be $[\hd(M)]$.
	}
\end{notation}

\begin{proposition}\label{prop: preserving green corr and head}
	Keep the notation of \ref{notation: green corr and head}. For any $\sigma\in \Gal(k/\F_p)$ and $\varphi\in \Aut(G)$, we have 
	$${}^\sigma f_P([V])=f_P({}^\sigma[V]),~~~{}^\sigma g_P([U])=g_{P}({}^\sigma[U]),~~~{}^\sigma \hd([M])=\hd({}^\sigma[M]).$$
	$${}^\varphi f_P([V])=f_{\varphi(P)}({}^\varphi[V]),~~~{}^\varphi g_P([U])=g_{\varphi(P)}({}^\varphi[U]),~~~{}^\varphi \hd([M])=\hd({}^\varphi[M]).$$
\end{proposition}

\begin{proof}
	It is easy to see that the functors $\sigma$ and $\varphi$ are compatible with restriction, induction, taking a vertex, and taking radical. We also note that these two functors preserve the order of vertices of indecomposable modules. Hence by the definitions of a Green correspondent and a head, the statement holds.
\end{proof}

\section{On blocks with an abelian defect group and proof of Theorem \ref{theo:cyclic}}\label{s3}

%In this section, all modules considered are finite dimensional. 

\begin{lemma}\label{lem:dimension at least}
Let $G$ be a finite group, $b$ a block of $kG$ with a defect group $D$, and $c$ the Brauer correspondent of $b$ in $kN_G(D)$. For any nonzero $kGb$-module $V$, the $k N_G(D)c$-module $cV$ is nonzero. Consequently, if every simple $kN_G(D)c$-module has dimension $m$, then every simple $kGb$-module has dimension at least $m$.
\end{lemma}

\begin{proof}
By the definition of a Brauer correspondent (\cite[Theorem 6.7.1]{Lin18b}), $c=\br_D^{kG}(b)$. Consider the idempotent $bc\in (kG)^D$. Since $\br_D^{kG}(bc)=c\neq 0$, by \cite[Theorem 6.4.6]{Lin18b}, sending any $kGb$-module $M$ to the $bckGbc$-module $bcM$ induces a Morita equivalent between $kGb$ and $bckGbc$. Thus $cV=bcV\neq0$ as $V\neq 0$.

Suppose that every simple $k N_G(P)c$-module has dimension $m$. For any simple $kGb$-module $V$, $cV$ is nonzero by the preceding paragraph. Taking a simple composition factor of $cV$ yields $m \leq \dim_k(cV) \leq \dim_k(V)$.
\end{proof}

\begin{lemma}\label{lem:dimension}
Let $G$ be a finite group, $b$ a block of $kG$ with a defect group $D$, and $c$ the Brauer correspondent of $b$ in $kN_G(D)$. Let $e$ be a block of $kC_G(D)$ such that $\br_D^{kG}(b)e\neq 0$. Then the block $e$ has a central defect group $Z(D)$ (see \cite[Corollary 6.3.10]{Lin18b}), and there is a unique simple $kC_G(D)e$-module (say $T$) up to isomorphism (see \cite[Proposition 6.6.5]{Lin18b}). Set $E=N_G(D,e)/DC_G(D)$, $d =\dim_k(T)$ and $t=|E|$. Suppose that $D$ is abelian and $E$ is cyclic. Then $kN_G(D)c$ has exactly $t$ isomorphism classes of simple modules, all of dimension $|N_G(D):N_G(D,e)|d$. In particular, $$\dim_k([U]) = |N_G(D):N_G(D,e)|d \leq \dim_k([V]),$$
for any $[V] \in \IBr(G,b)$ and any $[U] \in \IBr(N_G(D),c)$.
\end{lemma}

\begin{proof}
By \cite[Theorem 6.2.6 (iii)]{Lin18b}, $e$ remains a block of $kN_G(D,e)$. Since $T$ is the unique simple $k C_G(P) e$-module up to isomorphism, its isomorphism class is stable under $N_G(P,e)$-conjugation. Since $E=N_G(D,e)/C_G(D)$ is cyclic, by \cite[Theorem 8.12]{Navarro}, $T$ extends to a simple $kN_G(P,e)$-module $\widehat{T}$.  By the Clifford theory we see that for any simple $kN_G(D,e)e$-module $V$, $T$ is isomorphic to a composition factor of $\Res_{C_G(D)}^{N_G(D,e)}(V)$. Conversely, if $V$ is a simple $kN_G(D,e)$-module such that $T$ is isomorphic to a composition factor of $\Res_{C_G(D)}^{N_G(D,e)}(V)$, then $V$ belongs to the block $e$ of $kN_G(D,e)$. Therefore, by \cite[Corollary 8.20]{Navarro}, we have $|\IBr(N_G(D,e),e)|=|\IBr(E)|$. Since $E$ is abelian and $p\nmid |E|$ (see e.g. \cite[Theorem 6.7.6 (v)]{Lin18b}), we have $|\IBr(E)|=|E|=t$. Then by \cite[Theorem 6.7.6 (iii)]{Lin18b}, we have $|\IBr(N_G(D),c)|=|\IBr(N_G(D,e),e)|=t$. Again by \cite[Corollary 8.20]{Navarro}, any simple $k N_G(P,e)e$-module is isomorphic to $\widehat{T} \otimes_k S$ for some simple $kE$-module $S$, where we regard $S$ as a simple $kN_G(D,e)$-module via the canonical homomorphism $N_G(P,e)\to E$. Since $E$ is abelian, for any simple $kE$-module we have $\dim_k(S)=1$. It follows by \cite[Theorem 6.7.6 (iii)]{Lin18b} that any simple $k N_G(P)c$-module is isomorphic to $\Ind_{N_G(P,e)}^{N_G(P)}(\widehat{T} \otimes_k S)$, for some simple $kE$-module $S$. Hence any simple $kN_G(P)c$-module is of dimension $|N_G(P):N_G(D,e)|d$. The last statement follows from Lemma \ref{lem:dimension at least}.
\end{proof}

\begin{proposition}[A refinement of {\cite[Consequences 1, 2]{Alp}}]\label{lem:c and weight}

Let $G$ be a finite group, $b$ a block of $kG$ with an abelian defect group $D$, and $c$ the Brauer correspondent of $b$ in $kN_G(D)$. Then there is a canonical bijection 
$$\bigsqcup_{P\in [s_p(G)/\sim_G]}\W(N_G(P)/P,\bar{b}_P)\to \IBr(N_G(D),c)$$ 
preserves dimensions of modules and commutes with the action of $\H_b$, where $\H=\Gal(k/\F_p)\times \Aut(G)$.
\end{proposition}	

\begin{proof}
A slightly weaker version of this proposition, which only considers Galois automorphisms, has been used in \cite[page 324]{HY}. By modifying the choice of $[s_p(G)/\sim_G]$, we may assume that $D\in [s_p(G)/\sim_G]$. Hence $D$ is abelian. Then by (the proof of) \cite[Consequence 2]{Alp}, we have
$$\bigsqcup_{P\in [s_p(G)/\sim_G]}\W(N_G(P)/P,\bar{b}_P)=\W(N_G(D)/D,\bar{b}_D).$$
By (the proof of) \cite[Consequence 1]{Alp}, the inflation map $\W(N_G(D)/D,\bar{b}_D)\to \IBr(N_G(D),c)$ is a bijection. Using the definition of the inflation map, this bijection preserves dimensions of modules and commutes with the action of $\H_b$.
\end{proof}

%\section{Proof of Theorems \ref{theo:cyclic} and \ref{theo:Klein four}}

\begin{proof}[Proof of Theorem \ref{theo:cyclic} for cyclic defect groups]
By Proposition \ref{lem:c and weight}, it suffices to show that there is a bijection
$$f:\IBr(G,b)\to\IBr(N_G(D),c)$$
satisfying the required two conditions.

Since $D$ is cyclic, by \cite[Lemma~4.1]{KS}, there is a bijection
$$f:\IBr(G,b)\longrightarrow\IBr(N_G(D),c), \qquad [V] \longmapsto \hd(f_D([V])).$$
By \cite[Theorem 11.1.1]{Lin18b}, the block $b$ satisfying the assumption of Lemma \ref{lem:dimension}. By Lemma \ref{lem:dimension}, for any $[U]\in \IBr(N_G(D),c)$ and any $[V]\in \IBr(G,b)$, we have $\dim_k([U])\leq \dim_k([V])$. Hence all bijections between $\IBr(G,b)$ and $\IBr(N_G(D),c)$, including $f$, satisfy the required condition (i).
By Proposition \ref{prop: preserving green corr and head}, the bijection $f$ also satisfies the required condition (ii), whence the result. 
\end{proof}

Let $G$ be a finite group, $b$ a block of $kG$ with a defect group $D$. Let $e$ be a block of $kC_G(D)$ such that $\br_D^{kG}(b)e\neq 0$. The group $N_G(D,e)/DC_G(D)$ is a $p'$-group (see e.g. \cite[Theorem 6.7.6 (v)]{Lin18b}), and it is called a {\it inertial quotient} of $b$. Following Puig \cite[Definition 3.2]{Puig81}, a primitive idempotent $i\in(kGb)^D$ satisfying $\br_D^{kG}(i)\neq 0$ is called a {\it source idempotent} of $b$, and the $k$-algebra $ikGi$ is called a {\it source algebra} of $b$. By a result of Puig (see \cite[Corollary 3.5]{Puig81} or \cite[Theorem 6.4.6]{Lin18b}), there is canonical Morita equivalence between $kGb$ and $ikGi$, sending any $kGb$-module $V$ to the $ikGi$-module $iV$. %Following Kessar and Linckelmann \cite[Definition 1.1]{KL2025}, the $kGb$-module
%$kGbi \otimes_{k P} k$ is called a {\it source permutation module} of $kGb$.

\begin{proof}[Proof of Theorem \ref{theo:cyclic} for Klein four defect group]
By Proposition \ref{lem:c and weight}, it suffices to show that there is a bijection
$$f:\IBr(G,b)\to\IBr(N_G(D),c)$$
satisfying the required two conditions.

By \cite[Theorem 1.3]{KL2025}, for any simple $kN_G(D)c$-module $U$, its Green correspondent $g_D(U)$ is isomorphic to a direct summand of $kGi\otimes_{kD}k$. In particular, $kGi\otimes_{kD}k$ has at least $|\IBr(N_G(D),c)|$ isomorphism classes of indecomposable direct summands.
By \cite[Theorem 1.4]{KL2025}, any simple $kGb$-module occurs both as a quotient and as a submodule of $kGi\otimes_{kD} k$, and hence all simple $kGb$-modules are contained in both the head and the socle of $kGi\otimes_{kD}k$.

Since $D$ is a Klein four group, by \cite[Theorem 1.10]{KL2025}, the $k$-algebra $\End_{kGb}(kGi\otimes_{kD} k)$ is self-injective.
By the ``moreover" part of \cite[Theorem 2.13 (b)]{Geck}, taking heads gives a bijection between the isomorphism classes of indecomposable direct summands of $kGi\otimes_{kD}k$ and $\IBr(G,b)$. By Erdmann's theorem in \cite{Erdmann}, any simple module of a block with a Klein four defect group $D$ has $D$ as a vertex. Thus, there is a well defined map
$$\Omega:\IBr(N_G(D),c) \longrightarrow \IBr(G,b),\qquad [U]\longmapsto \hd(g_D([U])),$$
and $\Omega$ is injective. It is known that $|\IBr(G,b)|=|\IBr(N_G(D),c)|$ (Alperin's weight conjecture was known for Klein four defect groups). Hence $\Omega$ is a bijection. 

Let $f$ be the inverse of $\Omega$. Since any inertial quotient $E$ of $b$ is a $p'$-subgroup of $\Aut(D)$, $E$ is cyclic. By Lemma \ref{lem:dimension}, $f$ satisfies the required condition (i). By Proposition \ref{prop: preserving green corr and head}, $\Omega$ commutes with the action of $\H_b$. Thus $f$ also satisfies the required condition (ii), whence the result.
\end{proof}

\section{On inertial blocks and proof of Theorem \ref{theo:inertial blocks}}\label{s4}

\begin{proposition}\label{prop:inertial block, c and weight}
Let $G$ be a finite group, $b$ a block of $kG$ with a defect group $D$, and $c$ the Brauer correspondent of $b$ in $kN_G(D)$. Assume that there is a Morita equivalence between $kGb$ and $kN_G(D)c$ induced by a bimodule with an endopermutation source. Then there is a canonical bijection 
$$\bigsqcup_{P\in [s_p(G)/\sim_G]}\W(N_G(P)/P,\bar{b}_P)\to \IBr(N_G(D),c)$$ 
preserves dimensions of modules and commutes with the action of $\H_b$, where $\H=\Gal(k/\F_p)\times \Aut(G)$.
\end{proposition}

\begin{proof}
By modifying the choice of $[s_p(G)/\sim_G]$ we may assume that $D\in [s_p(G)/\sim_G]$. By (the proof of) \cite[Consequence 1]{Alp} the inflation map $\W(N_G(D)/D,\bar{b}_D)\to \IBr(N_G(D),c)$ is a bijection. Using the definition of the inflation map, this bijection preserves dimensions of modules and commutes with the action of $\H_b$.
Hence to prove the statement, it suffices to show that 
$$\bigsqcup_{P\in [s_p(G)/\sim_G]}\W(N_G(P)/P,\bar{b}_P)=\W(N_G(D)/D,\bar{b}_D).$$
For this, by \cite[Corollary 1.3]{H25}, we may assume that $G=N_G(D)$. Assume that there is $P\in [s_p(G)/\sim_G]$ such that $\W(N_G(P)/P,\bar{b}_P)\neq\varnothing$. Then $b_P=\br_P^{kG}(b)\neq 0$, and hence $P$ is conjugate to a subgroup of $D$. By modified the choice of $[s_p(G)/\sim_G]$ we may assume that $P\leq D$. Assume that $P<D$. Since $kN_G(P)/P$ has a projective simple module, by \cite[Proposition 1.15.14]{Lin18a}, we have $O_p(N_G(P))=P$. Since $D\unlhd G$, we have $N_D(P)=D\cap N_G(P)\unlhd N_G(P)$, and hence $N_D(P)\leq O_p(N_G(P))=P$. This contradicts to a well-known group-theoretic fact that $N_D(P)>P$. Now the proof is complete.
\end{proof}

\begin{lemma}[A version of {\cite[Lemma 3.2]{Lin25}} over positive characteristic]\label{lemma:inequality and source algebra}
Let $G$ be a finite group, $b$ a block of $kG$ with a defect group $D$, and $c$ the Brauer correspondent of $b$ in $kN_G(D)$. Let $i\in (kGb)^D$ and $j\in (kN_G(D)c)^D$ be source idempotents of $b$ and $c$, respectively. Let $X$ be a simple $kGb$-module and $Y$ a simple $kN_G(D)c$-module. If $\dim_k(jY)\leq \dim_k(iX)$, then $\dim_k(Y)\leq \dim_k(cX)\leq \dim_k(X)$.
\end{lemma}

\begin{proof}
The proof is similar to the proof of \cite[Lemma 3.2]{Lin25}.  By the definition of a Brauer correspondent (\cite[Theorem 6.7.1]{Lin18b}) and \cite[Theorem 6.2.6 (ii)]{Lin18b}, we see that $c=\br_D^{kN_G(D)}(c)$. Hence $c$ is a central idempotent in $kC_G(D)$. By \cite[Theorem 6.14.1]{Lin18b}, any primitive idempotent in $kC_G(D)c$ is a source idempotent in $(kN_G(D)c)^D$. Therefore, if $J$ is a primitive decomposition of $c$ in $kC_G(D)c$, then $J$ is also a primitive decomposaition of $c$ in $(kN_G(D)c)^D$, and every idempotent in $J$ is a source idempotent of $kN_G(D)c$. It follows that
$$\dim_k(Y)=|J|\cdot\dim_k(jY)~~~{\rm and}~~~\dim_k(cX)=|J|\cdot \dim_k(jX).$$ 
Note that $jb$ is an idempotent in $(kGb)^D$ such that $$\br_D^{kG}(jb)=\br_D^{kG}(j)\br_D^{kG}(b)=\br_D^{kG}(j)c=\br_D^{kG}(j)=\br_D^{kN_G(D)}(j)\neq 0.$$ 
Since $j$ is a primitive idempotent in $(kN_G(D))^D$ and $\br_D^{kN_G(D)}:(kN_G(D))^D\to k C_{N_G(D)}(D)= k C_G(D)$ is a surjective algebra homomorphism, $\br_D^{kN_G(D)}(j)$ is primitive in $kC_G(D)$; see e.g. \cite[Theorem 4.7.1 (ii)]{Lin18a}. Hence $bj=i'+i''$ for some source idempotent $i'\in (kGb)^D$ and some $i''\in (kGb)^D$ which is orthogonal to $i'$ and satisfies $\br_D^{kG}(i'')=0$. It follows that 
$$\dim_k(jX)=\dim_k(jbX)=\dim_k(i'X)+\dim_k(i''X)\geq \dim_k(i'X)=\dim_k(iX).$$
Thus we have
$$\dim_k(Y)=|J|\cdot\dim_k(jY)\leq |J|\cdot\dim_k(iX)\leq |J|\cdot\dim_k(jX)=\dim_k(cX),$$
as required.
\end{proof}

Let $G$ be a finite group. Recall that an {\it interior $G$-algebra} $A$ over $k$ is a $k$-algebra $A$ endowed with a group homomorphism $\rho:G\to A^\times$. This notion was defined by Puig \cite[Definition 3.1]{Puig81}. For example, if $M$ is a $kG$-module, then $\End_k(M)$ is an interior $G$-algebra via the structural homomorphism $G\to \End_k(M)^\times$. The group algebra $kG$ is an interior $G$-algebra in an obvious way. For any interior $G$-algebra $A$ and $H$ a subgroup of $G$, we set $A^H:=\{a\in A\mid \rho(h)a=a\rho(h)~{\rm for~all}~h\in H\}$; if $i\in A^H$ is an idempotent, then $iAi$ is an interior $H$-algebra with the structure homomorphism $H \to (iAi)^\times$ sending any $h \in H$ to $\rho(h)i$. If $A$ and $B$ are two interior $G$-algebras over $k$, then $A\otimes_k B$ is an interior $G$-algebra with the structural homomorphism
$$G\to G\times G\to A^\times\times B^\times\to (A\otimes_k B)^\times,$$
where the first map $G\to G\times G$ sends any $g\in G$ to $(g,g)\in G\times G$.
A {\it homomorphism} of two $G$-interior algebras is a $k$-algebra homomorphism compatible with the structure of interior $G$-algebra. For a homomorphism $f:G\to H$ of finite groups and an interior $H$-algebra $A$, we denote by ${}_fA$ the interior $G$-algebra which is equal to $A$ as a $k$-algebra, endowed with the structural homomorphism $G\xrightarrow{f} H\xrightarrow{\rho} A^\times$.

\begin{proposition}[cf. {\cite[2.16]{Puig11}}]\label{prop: source algebra of inertial block}
Let $G$ be a finite group, $b$ a block of $kG$ with a defect group $D$, $c$ the Brauer correspondent of $b$ in $kN_G(D)$. Let $i\in (kGb)^D$ and $j\in (kN_G(D)c)^D$ be source idempotents. Assume that there is a $kGb$-$kN_G(D)c$-bimodule $M$ having an endopermutation module as a source and inducing a Morita equivalence between $kGb$ and $kN_G(D)c$. Then there is an indecomposable endopermutation $kD$-module $W$ such that 
$$ikGi\cong \End_k(W)\otimes_k jkN_G(D)j$$
as $k$-algebras.
\end{proposition}

\begin{proof}
Write $A=ikGi$ and $B=jkN_G(D)j$. By the canonical Morita equivalence between a block algebra and its source algebra, the $A$-$B$-bimodule $iMj$ induces a Morita equivalence between $A$ and $B$. Then by standard Morita theory, the structural homomorphism $A\to\End_{B^{op}}(iMj)$ is an isomorphism of $k$-algebras, where $B^{\rm op}$ denotes the opposite $k$-algebra of $B$. This $k$-algebra isomorphism is automatically an interior $D$-algebra isomorphism. Since $i$ is a primitive idempotent in $A^D$, the identity map of $iMj$ is a primitive idempotent in $\End_{kD\otimes_k B^{\rm op}}(iMj)$. By \cite[Theorem 9.11.2]{Lin18b}, there is an automorphism $\varphi\in \Aut(D)$ such that $M$ has a vertex $\Delta \varphi:=\{(u,\varphi(u))\mid u\in D\}$. The group $\Delta \varphi$ is called a {\it twisted diagonal subgroup} of $G\times N_G(D)$.
Then by \cite[Proposition 3.4 (i)]{H26a}, there is an indecomposable endopermutation $kD$-module $W$, a group isomorphism $\tau:D\to D$, and some primitive idempotent $e\in (\End_k(W)\otimes_k {}_\tau B)^D$ such that 
$$A\cong e(\End_k(W)\otimes_k {}_\tau B)e$$
as interior $D$-algebras. By \cite[Theorem A]{Kul} or \cite[Theorem 6.14.1]{Lin18b}, $B\cong k_\alpha (D\rtimes E)$ as interior $D$-algebras, where $E$ is a $p'$-subgroup of $\Aut(D)$, and where $\alpha\in H^2(E;k^\times)$, inflated to $D\rtimes E$ and regarded as an element of $H^2(D\rtimes E;k^\times)$. Then by \cite[Corollary 7.4.5]{Lin18b} we see that the unity element of $(\End_k(W)\otimes_k {}_\tau B)^D$ is a primitive idempotent, and hence $e=1$. This shows that 
$$A\cong \End_k(W)\otimes_k {}_\tau B$$
as interior $D$-algebras. Since ${}_\tau B\cong B$ as $k$-algebras, forgetting the structure of interior $D$-algebras we obtain
$$A\cong \End_k(W)\otimes_k B$$
as $k$-algebras.
\end{proof}

\begin{proof}[Proof of Theorem \ref{theo:inertial blocks}]
By Proposition \ref{prop:inertial block, c and weight}, it suffices to show that there is a bijection 
$$f:\IBr(G,b)\to\IBr(N_G(D),c)$$ 
such that for any $[V]\in \IBr(G,b)$, we have $\dim_k(f([V]))\leq \dim_k([V])$. Let $i\in(kGb)^D$ be a source idempotent of $b$.
By Proposition \ref{prop: source algebra of inertial block}, the source algebra $ikGi$ is isomorphic to $\End_k(W)\otimes_k jkN_G(D)j$ as $k$-algebras, where $W$ is an endopermutation $kD$-module and $j\in (kN_G(D)c)^D$ is a source idempotent of $c$. For any $jk N_G(D) j$-module $T$, we can regard $W\otimes_kT$ as a $ikGi$-module via the above isomorphism; the attachment $T\mapsto W\otimes_k T$ induces a Morita equivalence between $jk N_G(D) j$ and $ikGi$. By the canonical Morita equivalence between a block algebra and its source algebra, for any simple $kGb$-module $V$, there is, up to isomorphism, a unique simple $kN_G(D)c$-module $U$ such that $iV\cong W\otimes_k jU$ as $ikGi$-modules. This defines a bijection 
$$f:\IBr(G,b)\to\IBr(N_G(D),c)$$ 
sending $[U]$ to $[V]$.
Since $\dim_k(iV)=\dim_k(W)\dim_k(jU)\geq \dim_k(jU)$, by Lemma \ref{lemma:inequality and source algebra}, we have $\dim_k(U)\leq \dim_k(V)$. Hence $f$ satisfies the required condition.
\end{proof}

\section{On $p$-solvable groups and proof of Theorem \ref{theo:p-solvable groups}}\label{s5}

In this section, we prove Theorem \ref{theo:p-solvable groups} by showing that the bijections constructed in \cite[Theorem 7.1]{Turull14} can be chosen to satisfy the required dimension inequality.
We first recall some notation from \cite[Section 5]{Turull14}.

\begin{notation}
{\rm
Let $G$ be a finite group and let $Q$ be a $p$-subgroup of $G$. Denote by $\IBr(G,Q)$ the set of isomorphism classes of simple $kG$-modules having $Q$ as a vertex. Denote by $\W(G,Q)$ the set of isomorphism classes of simple $k N_G(Q)$-modules obtained by inflation from projective simple $k N_G(Q)/Q$-modules.

Let $U$ be a normal $p'$-subgroup of $G$. Since $U$ is a $p'$-group, there is a canonical bijection between $\IBr(U)$ and the set of irreducible ordinary characters of $U$ (see \cite[Theorem 2.12]{Navarro}). Consider the conjugation action of $Q$ on $U$. Let 
$$\IBr(U)^Q:=\{[M]\in \IBr(U)\mid {}^{c_u}[M]=[M]~{\rm for~all}~u\in Q\}.$$ 
By \cite[Theorem 2.9]{Navarro2018}, there is a bijection 
$$\gamma_Q:\IBr(U)^Q\to \IBr(C_U(Q)),$$
which is called a {\it Glauberman correspondence}. Let $[M]\in \IBr(U)^Q$ and let $M_Q$ be a simple $kC_U(Q)$-module such that $[M_Q] = \gamma_Q ([M])$.
For any subgroup $H$ of $G$ containing $U$, we denote by $\IBr(H\mid M)$ the set of isomorphism classes of simple $kH$-modules whose restriction to $U$ has a composition factor (which is in fact a direct summand) isomorphic to $M$.
Define the {\it inertia group} $I_G(M)$ of $M$ in $G$ as
$$I_G(M) = \{g\in G \mid {}^{c_g}[M]=[M]\}$$
and the {\it extended inertia group} $\widetilde{I}_G(M)$ of $M$ in $G$ as 
$$\widetilde{I}_G(M) = \{g\in G\mid {}^{c_g}[M]={}^{\sigma}[M]
\text{ for some }\sigma\in\Gal(k/\F_p)\}.$$
Then $Q\leq I_G(M) \unlhd \widetilde{I}_G(M)$. Following \cite[Definitions 5.5 and 5.6]{Turull14}, we define
$$\IBr(G,Q\mid M,\F_p) = \{ [\Ind_{I_G(M)}^G(V)] \mid [V]\in\IBr(I_G(M) \mid {}^{\sigma} M)~\text{for some }$$ 
$$\sigma\in \Gal(k/\F_p)~\text{such that}~V~\text{has a vertex}~Q\}$$
and
\[
\W(G,Q\mid M_Q,\F_p) = \W(G,Q) \cap \IBr(N_G(Q),Q \mid M_Q , \mathbb{F}_p).
\]
By the equivalence of the Galois actions described in \cite[\S 2 (2")]{Turull14}, these two sets are canonically corresponding to the sets defined in \cite[Definitions 5.5 and 5.6]{Turull14}, respectively. Therefore, by \cite[Lemma 4.1]{Turull14}, $\IBr(G,Q\mid M,\F_p)$ is a subset of $\IBr(G,Q)$.
}
\end{notation}
The following lemma provides the essential dimension inequality in our proof.

\begin{lemma}\label{lem:psolvable-turull}
Let $G$ be a finite $p$-solvable group and $Q$ a $p$-subgroup of $G$. Let $U=O_{p'}(G)$ and $[M] \in \IBr(U)^Q$. Suppose that $\widetilde{I}_G(M)=G$. Set $I = I_G(M)$ and $R=Q\cap O_{p',p}(I)$. Suppose that $R\in\mathrm{Syl}_p(O_{p',p}(I))$. Let $M_R$ be a simple $kC_U(R)$-module such that $[M_R]=\gamma_R([M])$. Then there exists a bijection
\[
\tau: \IBr(G,Q\mid M,\F_p) \longrightarrow \IBr(N_G(R),Q\mid M_R,\F_p)
\]
satisfying the following conditions:
\begin{enumerate}[{\rm (i)}]
\item For any $[V]\in\IBr(G,Q\mid M,\F_p)$, we have $\dim_k(\tau([V])) \leq \dim_k([V])$.

\item For any $[V]\in\IBr(G,Q\mid M,\F_p)$, we have ${\rm codeg}([V])_{p'}^2 \equiv {\rm codeg}(\tau([V]))_{p'}^2 \pmod{p}$.

\item Let $b$ be a block of $kG$ and $[V]\in\IBr(G,b)\cap\IBr(G,Q\mid M,\F_p)$. If $\tau([V])=[S]$, then $b_R S=S$; recall that $b_R:=\br_R^{kG}(b)$.

\item The bijection $\tau$ is commutative with the action of $\Gal(k/\F_p)$.
\end{enumerate}
\end{lemma}

\begin{proof}
By \cite[Theorem 6.1 and \S 2(2")]{Turull14}, there is a bijection 
\[
\tau: \IBr(G,Q\mid M,\F_p) \longrightarrow \IBr(N_G(R),Q\mid M_R,\F_p)
\]
and $\tau$ is commutative with the action of $\Gal(k/\F_p)$ and satisfies condition (ii). Set $q= \dim_k(M_R) / \dim_k(M)$. For any $[V] \in \IBr(G,Q\mid M,\F_p)$, set $\tau([V])=[S]$. By the proof of \cite[Theorem 6.1]{Turull14}, there exists an equality $\dim_k(S) = q \dim_k(V)$. By the Glauberman correspondence, $M_R$ occurs as a composition factor of $\Res_{C_U(R)}^U(M)$. Hence $0<q \leq 1$. Consequently, $\dim_k(\tau([V]))\leq\dim_k([V])$.
Suppose that $[V]\in\IBr(G,b)$ and $[S]\in\IBr(N_G(R),c)$, where $b$ and $c$ are blocks of $kG$ and $kN_G(R)$, respectively.
By \cite[Theorem 6.1(1)]{Turull14}, we have $\br_R^{kG}(b)c=c$, and hence $b_RS=S$, as required.
\end{proof}

\begin{theorem}[A refinement of {\cite[Theorem 7.1]{Turull14}}]\label{theo:psolvable-relative}
Let $G$ be a finite $p$-solvable group and $Q$ a $p$-subgroup of $G$. Let $U$ be a normal $p'$-subgroup of $G$ and $[M] \in \IBr(U)^Q$. Let $M_Q$ be a simple $kC_U(Q)$-module such that $[M_Q] = \gamma_Q([M])$. Then there exists a bijection
\[
f:\IBr(G,Q\mid M,\F_p)\longrightarrow\W(G,Q\mid M_Q,\F_p)
\]
satisfying the following conditions:
\begin{enumerate}[{\rm (i)}]
\item For any $[V]\in\IBr(G,Q\mid M,\F_p)$, we have $\dim_k(f([V]))\leq\dim_k([V])$.

\item For any $[V]\in\IBr(G,Q\mid M,\F_p)$, we have ${\rm codeg}([V])_{p'}^2 \equiv {\rm codeg}(f([V]))_{p'}^2 \pmod{p}$.

\item Let $b$ be a block of $kG$ and $[V]\in\IBr(G,b)\cap\IBr(G,Q\mid M,\F_p)$. If $f([V])=[S]$, then $b _Q S = S$; recall that $b_Q := \br_Q^{kG}(b)$.

\item The bijection $f$ is commutative with the action of $\Gal(k/\F_p)$.
\end{enumerate}
\end{theorem}

\begin{proof}
We follow the inductive construction in the proof of \cite[Theorem 7.1]{Turull14}, adding the dimension inequality condition. By \cite[Theorem 7.1 and \S 2(2")]{Turull14}, the bijections constructed there commute with the action of $\Gal(k/\F_p)$ and satisfy condition (ii). Moreover, if $[V]\in\IBr(G,b)$ and $f([V])=[S]\in\IBr(N_G(Q),c)$, then $\br_Q^{kG}(b)c=c$ by \cite[Theorem 7.1(1)]{Turull14}. Therefore $b_QS=S$. Thus it remains to show that the bijections in the inductive construction can be chosen to satisfy the required dimension inequality.

We argue by induction on $(|G|,|G:U|)$, ordered lexicographically. Suppose that the assertion is false, and choose a counterexample minimal with respect to this ordering. Suppose that $Q=1$. By \cite[Proposition 5.8]{Turull14}, $\W(G,1\mid M_1,\F_p) = \IBr(G,1 \mid M,\F_p)$ and thus the identity map gives the required bijection. Hence $Q\neq1$.

Suppose that the extended inertia group $\widetilde{I}_G(M) < G$. By \cite[Propositions 5.9 and 5.10]{Turull14}, inductions yield bijections
\[
b_1 = \Ind_{\widetilde{I}_G(M)}^G:\IBr(\widetilde{I}_G(M),Q\mid M,\F_p) \longrightarrow \IBr(G,Q\mid M,\F_p),
\]
and
\[
b_2 = \Ind_{N_{\widetilde{I}_G(M)}(Q)}^{N_G(Q)}:\W(\widetilde{I}_G(M),Q\mid M_Q,\F_p)
\longrightarrow
\W(G,Q\mid M_Q,\F_p).
\]
By induction hypothesis, there exists a bijection
\[
f_{\widetilde{I}_G(M)}:\IBr(\widetilde{I}_G(M),Q\mid M,\F_p)
\longrightarrow
\W(\widetilde{I}_G(M),Q\mid M_Q,\F_p)
\]
satisfying all the required properties. Define $f=b_2\circ f_{\widetilde{I}_G(M)}\circ b_1^{-1}$. Let $[V]\in\IBr(\widetilde{I}_G(M),Q \mid M,\F_p)$ and set $[S]=f_{\widetilde{I}_G(M)}([V])$. By induction hypothesis, $\dim_k(S)\leq\dim_k(V)$. This gives
\[
|N_G(Q):N_{\widetilde{I}_G(M)}(Q)|\dim_k(S)
\leq
|N_G(Q):N_{\widetilde{I}_G(M)}(Q)|\dim_k(V)
\leq
|G:\widetilde{I}_G(M)|\dim_k(V).
\]
Thus $f$ satisfies the required dimension inequality, which is a contradiction. Hence $\widetilde{I}_G(M) = G$. In particular, $I_G(M)$ is normal in $G$.

Suppose that $U<O_{p'}(G)$. Set $W=O_{p'}(G)$. Since $\widetilde{I}_G(M) = G$, the group $N_G(Q) \times \Gal(k/\F_p)$ acts on $\IBr(C_W(Q), 1 \mid M_Q , \F_p)$. Let $O$ be the set of orbits of this action. For each $o \in O$, choose a representative $[N_{Q,o}]\in o$. Let $N_o$ be a simple $k W$-module such that $[N_{Q,o}] = \gamma_Q ([N_o])$. As in the second paragraph of the proof of \cite[Theorem 7.1]{Turull14}, we have disjoint unions
\[
\IBr(G,Q \mid M , \F_p) = \bigsqcup_{o \in O} \IBr(G,Q \mid N_o , \F_p)
\]
and
\[
\W(G,Q \mid M_Q , \F_p) = \bigsqcup_{o \in O} \W(G,Q \mid N_{Q,o}, \F_p).
\]
Since $|G:W|<|G:U|$, the induction hypothesis applies to each $N_o$. Thus, for each $o\in O$, there exists a bijection
$f_o: \IBr(G,Q \mid N_o,\F_p) \to \W(G,Q \mid N_{Q,o},\F_p)$ satisfying $\dim_k(f_o([V])) \leq\dim_k([V])$ for any $[V] \in \IBr(G,Q \mid N_o,\F_p)$. Taking the disjoint union of the bijections $f_o$ gives the required bijection $f$. This contradiction shows that $U=O_{p'}(G)$.

Set $I=I_G(M)$ and $R=Q\cap O_{p',p}(I)$. Suppose that $R$ is not a Sylow $p$-subgroup of $O_{p',p}(I)$. By the argument in the third paragraph of the proof of \cite[Theorem 7.1]{Turull14}, both $\IBr(G,Q \mid M,\F_p)$ and $\W(G,Q \mid M_Q,\F_p)$ are empty. The trivial bijection yields a contradiction. Thus $R$ is a Sylow $p$-subgroup of $O_{p',p}(I)$.

Suppose that $H=N_G(R)<G$. Let $M_R$ be a simple $kC_U(R)$-module such that $[M_R] = \gamma_R ([M])$. By Lemma \ref{lem:psolvable-turull}, there exists a bijection
\[
\tau: \IBr(G,Q\mid M,\F_p) \longrightarrow
\IBr(H,Q\mid M_R,\F_p)
\]
satisfying all the conditions of Lemma \ref{lem:psolvable-turull}. Since $I$ is normal in $\widetilde{I}_G(M)=G$ and $R = Q \cap O_{p',p}(I)$, we have $R \triangleleft N_G(Q)$ and hence $N_H(Q)=N_G(Q)$. Furthermore, by \cite[Theorem 2.9 and Problems (2.7)]{Navarro2018}, $\gamma_Q(\gamma_R([M])) = \gamma_Q([M]) = [M_Q]$ and $\W(H,Q\mid M_Q,\F_p) = \W(G,Q\mid M_Q,\F_p)$ by definition. Since $|H|<|G|$, the induction hypothesis provides a bijection $f_H: \IBr(H,Q\mid M_R,\F_p) \to \W(H,Q\mid M_Q,\F_p)$ satisfying the required conditions.

Define $f=f_H\circ\tau$. Then, for every $[V]\in\IBr(G,Q\mid M,\F_p)$, we have $\dim_k(f([V])) \leq \dim_k(\tau([V])) \leq\dim_k([V])$. Let $b$ be a block of $kG$ and suppose that $[V]\in\IBr(G,b)$. Write $\tau([V])=[T]\in\IBr(H,c)$ and $f_H([T])=[S]$, where $c$ is a block of $kH$. By Lemma \ref{lem:psolvable-turull}, we have $b_R T=T$. Since $b_R$ is a central idempotent of $kH$, it follows that $b_Rc=c$. Write $\br_Q^{kH} (c) = c_Q$. Using $\br_Q^{kH}(b_R)=\br_Q^{kG}(b)=b_Q$, we obtain $b_Qc_Q=c_Q$. By induction hypothesis, $c_Q S=S$, and hence $b_QS=S$.
By Lemma \ref{lem:psolvable-turull} and the induction hypothesis, both $\tau$ and $f_H$ satisfy conditions (ii) and (iv). Note that the conditions (ii) and (iv) are preserved under the composition. Thus $f$ satisfies all the required conditions, a contradiction. Therefore $R\triangleleft G$.

Suppose that $R\neq1$. Set $\overline G=G/R, \overline Q=Q/R, \overline U=UR/R$. Since $\overline U \cong U$, let $\overline M$ be the simple $k\overline U$-module corresponding to $M$ under this isomorphism. By \cite[Lemma 2.32]{Navarro}, inflation induces a bijection $\pi_1:\IBr(\overline G,\overline Q \mid \overline M,\F_p) \to \IBr(G,Q \mid M,\F_p)$. Since $C_{\overline U}(\overline Q)  = C_U(Q) R/R \cong C_U(Q)$, let $\overline M _{\overline Q}$ be the simple $k C_{\overline U}(\overline Q)$-module corresponding to $M_Q$. Then $\gamma_{\overline{Q}}([\overline M]) = [\overline M _{\overline Q}]$. Furthermore, $N_{\overline G}(\overline Q)/ \overline Q \cong N_G(Q)/Q$. Thus inflation also induces a bijection $\pi_2: \W(\overline G,\overline Q \mid \overline M_{\overline Q},\F_p) \to \W(G,Q \mid M_Q,\F_p)$.
By induction hypothesis, there exists a bijection $\overline f: \IBr(\overline G,\overline Q \mid \overline M,\F_p) \to \W(\overline G,\overline Q \mid \overline M_{\overline Q},\F_p)$.
Define $f=\pi_2\circ\overline f\circ\pi_1^{-1}$. Both bijections $\pi_1$ and $\pi_2$ preserve dimensions, commute with taking local blocks, and commute with the action of $\Gal(k/\F_p)$. Hence, for any $[V]\in\IBr(G,Q \mid M,\F_p)$, we have $\dim_k(f([V])) \leq\dim_k([V])$.
Then $f: \IBr(G,Q\mid M,\F_p) \longrightarrow \W(G,Q\mid M_Q,\F_p)$ is bijection satisfying all required conditions, a contradiction, as in the final reduction of the proof of \cite[Theorem 7.1]{Turull14}.

Finally, since $R$ is a Sylow $p$-subgroup of $O_{p',p}(I)$, the argument in the last paragraph of the proof of \cite[Theorem 7.1]{Turull14} yields $Q=1$, contrary to our assumption. This completes the induction and the proof.
\end{proof}

We are now ready to prove Theorem \ref{theo:p-solvable groups}.

\begin{proof}[Proof of Theorem \ref{theo:p-solvable groups}]
Let $G$ be a finite $p$-solvable group and let $b$ be a block of $kG$.

Taking $U=1$ and $M=k$ in Theorem \ref{theo:psolvable-relative}, for any $p$-subgroup $Q$ of $G$, we obtain a bijection
\[
f_Q:\IBr(G,Q) \longrightarrow \W(G,Q)
\]
satisfying all the conditions of Theorem \ref{theo:psolvable-relative}. Identifying $\W(G,Q)$ with $\W(N_G(Q)/Q)$ via inflation, and using $\IBr(G)=\bigsqcup_{Q\in[s_p(G)/\sim_G]}\IBr(G,Q)$, the bijections $f_Q$ induce a bijection
\[
f:\IBr(G)\longrightarrow \bigsqcup_{Q\in[s_p(G)/\sim_G]} \W(N_G(Q)/Q).
\]
Let $[V]\in\IBr(G,b)$ and suppose that $Q\in[s_p(G)/\sim_G]$ is a vertex of $V$. Write $f([V])=[S]$. By Theorem \ref{theo:psolvable-relative} (iii), we have $b_QS=S$. This is equivalent to $\bar b_QS=S$ as $k N_G(Q)/Q$-modules. Consequently, $[S]\in\W(N_G(Q)/Q,\bar b_Q)$. Conversely, let $[S]\in\W(N_G(Q)/Q,\bar b_Q)$ for some $Q\in[s_p(G)/\sim_G]$. Since $f$ is a bijection, there exists a unique $[V]\in\IBr(G,Q)$ such that $f([V])=[S]$. Suppose that $[V] \in \IBr(G,b')$ for a block $b'$ of $kG$. By Theorem \ref{theo:psolvable-relative} (iii), we have $b'_QS=S$. On the other hand, $b_QS=S$.
If $b\neq b'$, then $b_Qb'_Q =\br_Q^{kG}(bb')=0$ and $S=b_Qb'_QS=0$, a contradiction. Therefore $b=b'$ and $[V]\in\IBr(G,b)$.

Thus, $f$ restricts to a bijection
\[
f:\IBr(G,b)\longrightarrow \bigsqcup_{Q\in[s_p(G)/\sim_G]} \W(N_G(Q)/Q,\bar b_Q).
\]
Moreover, for every $[V]\in\IBr(G,b)$, Theorem \ref{theo:psolvable-relative} (i) and (ii) give $\dim_k(f([V]))\leq\dim_k([V])$ and ${\rm codeg}([V])_{p'}^2 \equiv {\rm codeg}(f([V]))_{p'}^2 \pmod{p}$. Since each $f_Q$ commutes with the action of $\Gal(k/\F_p)$ by Theorem \ref{theo:psolvable-relative} (iv), so does $f$. Consequently, its restriction to $\IBr(G,b)$ is commutative with the action of $\Gal(k/\F_p)_b$. This completes the proof.
\end{proof}

\appendix
\renewcommand{\thetheorem}{\Alph{section}.\arabic{theorem}}
\renewcommand{\theequation}{\Alph{section}.\arabic{equation}}

\section{}\label{representatives}

\begin{proposition}\label{prop:independence of representatives}
Let $G$ be a finite group. The 	$\Aut(G)$-set 
$$\bigsqcup_{P\in [s_p(G)/\sim_G]}\W(N_G(P)/P)$$ is independent, up to isomorphism, of the choice of the representatives $[s_p(G)/\sim_G]$.
\end{proposition}

\begin{proof}
Let $\mathcal{R}$ and $\mathcal{S}$ be two sets of representatives of the $G$-conjugacy classes of $p$-subgroups of $G$, and write
$$\W_{\mathcal{R}}=\bigsqcup_{P\in\mathcal{R}}\W(N_G(P)/P)~~~{\rm and}~~~\W_{\mathcal{S}}=\bigsqcup_{Q\in\mathcal{S}}\W(N_G(Q)/Q).$$
For each $P\in\mathcal{R}$, let $Q_P$ be the unique member of $\mathcal{S}$ conjugate to $P$, and choose $x_P\in G$ such that $Q_P={} ^{x_P}P$. Conjugation by $x_P$ induces an isomorphism $N_G(P)/P\to N_G(Q_P)/Q_P$. Transport of structure along this isomorphism preserves simplicity and projectivity, and hence defines a bijection
$$F:\W_{\mathcal{R}}\longrightarrow \W_{\mathcal{S}},\qquad [U]\in \W(N_G(P)/P)\longmapsto {}^{c_{x_P}}[U]\in \W(N_G(Q_P)/Q_P).$$
Any other choice of $x_P$ differs from it on the left by an element of $N_G(Q_P)$. Since inner automorphisms of $N_G(Q_P)/Q_P$ act trivially on isomorphism classes of modules, $F$ is independent of the choice of $x_P$.

Let $\varphi\in\Aut(G)$ and $[U]\in\W(N_G(P)/P)$, where $P\in\mathcal{R}$. Choose $P'\in\mathcal{R}$ and $g\in G$ such that $\varphi(P)={}^gP'$, and set $h=\varphi(x_P)g x_{P'}^{-1}$. Then $\varphi(Q_P)={}^hQ_{P'}$. By the definition of the actions described in \ref{void:The action of group automorphisms}~(v) and the identity $\varphi\circ c_{x_P}=c_{\varphi(x_P)}\circ\varphi$, we obtain
$$\begin{aligned}
{}^\varphi F([U])
&={} ^{c_h^{-1}\circ\varphi\circ c_{x_P}}[U]\\
&={} ^{c_{x_{P'}}\circ c_g^{-1}\circ\varphi}[U]\\
&=F({}^\varphi[U]).
\end{aligned}$$
Thus $F$ is an isomorphism of $\Aut(G)$-sets.
\end{proof}

\bigskip\noindent\textbf{Acknowledgements.} We appreciate Eugenio Giannelli for some very helpful comments, and for his refinement of the McKay conjecture, which inspired us to consider the simliar refinement for the Alperin weight conjecture.

The authors would like to thank ChatGPT for some helpful discussions. Specifically, we first observed that Conjecture \ref{conj:block version} holds for cyclic blocks, Klein four blocks, and inertial blocks, and then we asked GPT to find a conunterexample to Conjecture \ref{conj: non-block version}. No counterexample had been found. This encouraged us to further consider $p$-solvable groups. We consulted GPT for a few questions about $p$-solvable groups. These pieces of information were helpful, but were not directly used in this paper. The mathematical ideas, the proofs, and the writing of the paper were entirely developed by the authors. 

The authors acknowledge support from National Natural Scicence Foundation of China (12501024, 12471016) and China Postdoctoral Science Foundation (GZC20252006, 2025T001HB).

\bigskip
{\footnotesize School of Mathematics and Statistics, Central China Normal University, Wuhan 430079, China
	
Email address: xinhuang@mails.ccnu.edu.cn

\smallskip School of Mathematics and Statistics, Central China Normal University, Wuhan 430079, China

Email address: lihanxiao@mails.ccnu.edu.cn}

\end{document}